\documentclass[12pt]{amsart}
\usepackage{amsbsy}
\usepackage{graphicx,epsfig,subfigure,psfrag}
\begin{document}

\newtheorem{theorem}{Theorem}
\newtheorem{proposition}{Proposition}
\newtheorem{lemma}{Lemma}
\newtheorem{corollary}{Corollary}
\newtheorem{definition}{Definition}
\newtheorem{remark}{Remark}
\newcommand{\beq}{\begin{equation}}
\newcommand{\eeq}{\end{equation}}
\numberwithin{equation}{section}
\numberwithin{theorem}{section}
\numberwithin{proposition}{section}
\numberwithin{lemma}{section}
\numberwithin{corollary}{section}
\numberwithin{definition}{section}
\numberwithin{remark}{section}
\newcommand{\ren}{\mathbb{R}^N}
\newcommand{\re}{\mathbb{R}}
\newcommand{\n}{\nabla}
\newcommand{\iy}{\infty}
\newcommand{\pa}{\partial}
\newcommand{\fp}{\noindent}
\newcommand{\ms}{\medskip\vskip-.1cm}
\newcommand{\mpb}{\medskip}
\renewcommand{\a}{\alpha}
\renewcommand{\b}{\beta}
\newcommand{\g}{\gamma}
\newcommand{\G}{\Gamma}
\renewcommand{\d}{\delta}
\newcommand{\D}{\Delta}
\newcommand{\e}{\varepsilon}
\renewcommand{\l}{\lambda}
\renewcommand{\o}{\omega}
\renewcommand{\O}{\Omega}
\newcommand{\s}{\sigma}
\renewcommand{\t}{\tau}
\newcommand{\z}{z}
\newcommand{\wx}{\widetilde x}
\newcommand{\wt}{\widetilde t}
\newcommand{\noi}{\noindent}
\newcommand{\BB}{{\mathbf  B}}
\newcommand{\AAA}{{\mathbf  A}}
\newcommand{\CC}{{\mathbf  C}}
\newcommand{\ba}{\begin{eqnarray}}
\newcommand{\be}{\begin{equation}}
\newcommand{\ea}{\end{eqnarray}}
\newcommand{\ee}{\end{equation}}
\newcommand{\ssk}{\smallskip}
\def\com#1{\fbox{\parbox{6in}{\texttt{#1}}}}

\title
{\bf On nonexistence  of Baras--Goldstein type\\ for higher-order
 parabolic equations\\ with singular potentials}




\author {V.A.~Galaktionov and I.V.~Kamotski}

\address{Department of Mathematical Sciences, University of Bath,
 Bath BA2 7AY, UK}
\email{vag@maths.bath.ac.uk}

\address{Department of Mathematical Sciences, University of Bath,
 Bath BA2 7AY, UK}
\email{ivk20@maths.bath.ac.uk}


\keywords{Parabolic equations with singular potentials, Hardy
inequality, nonexistence, regular approximations, oscillatory
solutions.
$\qquad$ {\bf To appear in:} Trans. Amer. Math. Soc.}

\date{\today}
\begin{abstract}

The celebrated result by Baras and Goldstein (1984) established
that the heat equation with singular inverse square potential in a
smooth bounded domain $\O\subset \ren$, $N \ge 3$, such that $0
\in \O$,
 $$
  \mbox{$
 u_t= \D u + \frac c{|x|^2}\, u \quad \mbox{in}
 \quad \O \times (0,T), \quad u \big|_{\partial \O}=0,
  $}
  $$
in the supercritical range
 $$
  \mbox{$
 c> c_{\rm Hardy}(1)= \big(\frac {N-2}2\big)^2
 $}
 $$
does not have a solution for any nontrivial $L^1$  initial data
 $
 u_0(x) \ge 0$ in $\O
  $,
  or for a positive measure.
Namely,  it was proved that a regular approximation of a possible
solution by a sequence $\{u_n(x,t)\}$ of classical solutions of
uniformly parabolic equations with bounded truncated
 potentials given by
 $$
  \mbox{$
  V(x) = \frac c{|x|^2} \mapsto
 V_n(x)=\min \big\{ \frac c{|x|^2}, \, n \big\} \,\,\, (n \ge 1)
  $}
 $$
 diverges, and, as $n \to \infty$,
  $$
  u_n(x,t) \to +\infty
  \quad \mbox{in} \quad
  \O \times (0,T).
  $$

In the present paper, we reveal the connection of this ``very
singular" evolution with a spectrum of some ``limiting" operator.
The proposed approach allows us to consider more general
higher-order operators (for which Hardy's inequalities were known
since Rellich, 1954)  and initial data that are not necessarily
positive.
 In particular it is established that,
under some natural hypothesis, the divergence result is valid for
any $2m$th-order parabolic equation with singular potential
 $$
  \mbox{$
 u_t =-(-\D)^m u + \frac c{|x|^{2m}}\, u \quad \mbox{in}
 \quad \O \times (0,T), \quad \mbox{where} \quad c>c_{\rm H}(m),
 \,\,\, m \ge 1,
  $}
 $$
 with zero Dirichlet conditions on $\partial \O$ and for a wide class of
 initial data.
 In particular, typically, the
 divergence holds for any data satisfying
  $$
  u_0(x) \quad \mbox{is continuous at $x=0$\,\, and}
   \quad
  u_0(0)>0.
   $$

  Similar nonexistence (i.e., divergence as $\e \to 0$) results are also derived for time-dependent
  potentials $\e^{-2m}q(\frac x \e, \frac t{\e^{2m}})$ and
  nonlinear reaction terms $\frac{|u|^p}{\e^{2m}+|x|^{2m}}$ with $p>1$.
   Applications to other, linear and semilinear, Schr\"odinger and wave  PDEs are discussed.

\end{abstract}

\maketitle

\section{Introduction: Baras--Goldstein result  and  extensions}
   \label{Sect1}

The paper treats the questions of regular approximations of
higher-order parabolic equations and the corresponding elliptic
operators with singular unbounded potentials.
 There are two main generating and key ingredients
 of such a theory:

(i) The origin of such research can be attributed to Hardy's
inequalities (1919--20) for such symmetric operators in the
$L^2$-space, which established the parameter supercritical range,
where the operators are not semibounded and hence do not admit
Friedrichs' classic self-adjoint extensions (1935). It is also
crucial that, at the same time, in view of coinciding deficiency
indices, there exist infinitely many other self-adjoint extensions
with discrete spectra and $L^2$-eigenfunctions.

(ii) It was later the discovery of Baras--Goldstein (1984), that,
in the corresponding non-stationary parabolic flows, regular
approximations of nonnegative solutions uniformly {\em diverge},
i.e., the corresponding semigroup gives infinite values for any
such nontrivial data.

\ssk

It turned out later on that, in general singular linear and
nonlinear PDE theory, these are important problems concerning
existence of extended semigroups of proper (blow-up) solutions
that are obtained by regular approximations. For a number of
parabolic, hyperbolic, and other evolution equations of
mathematical physics, such semigroups can be essentially
discontinuous
 in any suitable or admissible metric. Actually, this means that
a proper solution blows up at some $t=T \ge 0$, ceases to exist as
a bounded solution for $t>T$, so a special framework of extended
semigroup theory via regularization should be put in charge
instead of the classic one.
 In particular,  proper setting of basic problems
such as the standard Cauchy or initial-boundary value ones
represents a difficult subject. Concerning related elliptic
operators with singular potentials, it is then key to understand,
which of their spectral properties admit a proper approximation by
a family of ``regularized" operators with truncated singularities.
It turns out that, even for classic Laplacian-type operators with
inverse square  potentials in the supercritical Hardy range,
standard self-adjoint extensions have nothing to do with the
actual operator that occurs in the regularization limit. It is
revealed that this actual ``limiting" operator has positive
eigenvalues and we demonstrate that they are responsible  for the
effective ``nonexistence" of the limit of a solution. This point
of view allows us not to rely on the Maximum Principle (which was
a standard tool in previous work related to Baras--Goldstein type
problems) and to consider much more general operators for which
the Maximum Principle is not available.  In particular, we
demonstrate it considering  higher-order elliptic operators. Such
$2m$th-order operators are  common in PDE theory and mathematical
physics, and the first derivation of Hardy's inequality for $m=2$
goes back to Rellich (1954).
 In what follows, we explain these aspects in
greater detail treating also other PDEs.

\ssk

 Thus, one of the  first
key results of modern theory of singular elliptic operators and
related extended discontinuous semigroups of blow-up solutions was
obtained  by Baras and Goldstein \cite{BG84}.

\subsection{Baras--Goldstein (1984): nonexistence for singular potential in linear heat equations}

Let $\O \subset \ren$ be a smooth bounded domain, where $N \ge 3$
and $0 \in \O$. In Baras--Goldstein  \cite{BG84}, the authors
considered
 the heat
equation with the inverse square potential  and the zero Dirichlet
boundary condition,
 \beq
 \label{e1}
  \mbox{$
 u_t= \D u + \frac c{|x|^2}\, u \quad \mbox{in}
 \quad \O \times (0,T), \quad u(x,t)=0 \,\,\, \mbox{on} \,\,\,
 \partial \O \times (0,T),
  $}
  \eeq
 and nonnegative initial
data
 \beq
 \label{e2}
 u(x,0)= u_0(x)\ge 0 \quad \mbox{in} \quad \O,
  \eeq
  where
  $u_0 \in L^1(\O)$, or, is a positive measure.

Their key {\em nonexistence} result is as follows:
  in the supercritical Hardy range
 \beq
 \label{e3}
  \mbox{$
 c> c_{\rm H}= \big(\frac {N-2}2\big)^2
 $}
 \eeq
the problem
 (\ref{e1}), (\ref{e2})
does not have a solution\footnote{The question on nonexistence was
posed to the authors by H.~Brezis and J.-L.~Lions,
\cite[p.~122]{BG84}.}.
  Namely, this means that the regular approximation of a possible solution by the sequence
$\{u_n(x,t)\}$ of bounded classical solutions satisfying
(\ref{e1}) with the truncated bounded potentials obtained on
replacement
 \beq
 \label{e5}
  \mbox{$
   V(x) = \frac c{|x|^2} \mapsto
 V_n(x)=\min \big\{ \frac c{|x|^2},\,\, n \big\}, \quad \mbox{with} \quad  n >0,
  $}
 \eeq
 diverges. More precisely, as $n \to \infty$,
   \beq
   \label{e6}
  u_n(x,t) \to +\infty \quad \mbox{in} \quad \O \times (0,T).
   \eeq
   On the contrary,
 for $c \le c_{\rm H}$ in (\ref{e1}), this sequence has a finite limit that
 corresponds to {\em existence} of a unique solution of the Cauchy problem
  (\ref{e1}), (\ref{e2}).


The results and ideas of the pioneering paper \cite{BG84}
generated a new direction of nonexistence/existence theory for
linear and nonlinear PDEs. We refer to \cite{GZ01, Gold02, Kom04}
for the study of linear parabolic equations and  to \cite{GGK05,
Kom06} devoted to quasilinear reaction-diffusion PDEs. These
questions are  reflected in the books \cite{MitPoh} and
\cite[Ch.~11]{GalGeom}.

\subsection{Main extensions to higher-order PDEs without
positivity assumptions on data}

As our main extended model, we consider the $2m$th-order parabolic
equation with a singular
potential,  
 \beq
 \label{2m1}
  \mbox{$
 u_t ={\bf B}_0 u \equiv -(-\D)^m u + \frac c{|x|^{2m}}\, u \quad \mbox{in}
 \quad \O \times \re_+, \quad \mbox{where} \quad
  m \ge 1, \,\,\, N>2m.
  $}
 \eeq
 For definiteness, we take
  zero Dirichlet conditions on the smooth boundary
 \beq
 \label{2m2}
  \mbox{$
 u= \frac{{\partial}u}{\partial \nu}=...=\frac{{\partial}^{m-1}u}{\partial
 \nu^{m-1}}=0 \quad \mbox{on} \quad \partial \O \times \re_+,
 $}
  \eeq
  where $\nu$ is the unit outward normal to $\partial \O$,
   and initial data
    \beq
    \label{2m3}
    u(x,0)=u_0(x) \in L^2(\O).
    \eeq
Actually, as one can expect, the nonexistence result is purely
associated with the strong singularity at $x=0$ and does not
essentially depend on boundary conditions (though their
self-adjoint nature is a convenient assumption for calculus
applied).

 As in Baras--Goldstein \cite{BG84}, our crucial
assumption is that the potential in (\ref{2m1}) belongs to  the
supercritical range, i.e., the constant $c$ is large enough:
 \beq
 \label{2m4}
 c>c_{\rm H}(m).
  \eeq
  Here $c_{\rm H}(m)$ for dimensions $N>2m$ is
 classic {\em Hardy's best constant} of multiplicative inequalities
involving the potential in (\ref{2m1}). This goes back to Hardy
(1919) for $m=1$, and Rellich (1954) for $m=2$;  see \cite{GGM}
and \cite{Yaf99} for further references and full history. The
Hardy constant is given by
 \beq
 \label{Best1N}
c_{\rm H}(m)=
 \left\{
 \begin{matrix}
B_2 B_4 ...B_m \quad \qquad\,\, \mbox{for $m$ even},
\smallskip\\
 B_3 B_5 ... B_m \, c_{\rm H}(1) \quad \mbox{for $m$ odd}.\,\,\,
 \end{matrix}
 \right.
 \eeq
 Here  $B_k = \big[\frac {(N-2k)(N+2k-4)}4\big]^2$ for $k=1,2,...,m$, and  $c_{\rm H}(1)$ is as
 in (\ref{e3}); see a simple derivation
 and a list of references in  \cite{GalH2}, and also \cite{BV, VZ} with a number of applications
  for $m=1$.

\smallskip

Similar to \cite{BG84} and as usual in extended semigroup theory
\cite{GV4, GalGeom}, we construct
 a {\em proper} solution of the problem (\ref{2m1})--({\ref{2m4}) using
regular approximations. For convenience, instead of (\ref{e5})
(clearly, this does not affect the final results), we perform an
analytic approximation of the potential by replacing
 \beq
 \label{a1}
  \mbox{$
V(x)= \frac c{|x|^{2m}} \mapsto  V_\e(x) = \frac c{\e^{2m} +
|x|^{2m}}, \quad \mbox{with} \quad \e>0.
  $}
  \eeq
  By $\{u_\e(x,t)\}$ we denote the  sequence of
  classical bounded (for $t>0$) solutions of this regularized initial-value
  problem for the parabolic PDE
  with the potentials (\ref{a1}),
 \beq
 \label{2m5}
  \mbox{$
  u_\e: \quad
 u_t ={\bf B}_\e u \equiv -(-\D)^m u + \frac c{\e^{2m}+|x|^{2m}}\, u \quad \mbox{in}
 \quad \O \times \re_+ \quad (m \ge 1, \,\,\,N>2m).
  $}
 \eeq
 We then take
  the same Dirichlet boundary conditions (\ref{2m2})
   and  same initial data
  (\ref{2m3}).

Passing to the limit as $\e \to 0$ is performed in Section
\ref{Sect01}, which is based on a more general approach in Section
\ref{SectG} to $2m$th-order parabolic equations with arbitrary
potentials. In particular, we show that the phenomenon:
   \beq
   \label{2m7}
    \{u_\e(x,t)\}\quad \mbox{diverges in $L^2(\O)$ as} \quad \e \to 0,
    \eeq
is a generic (robust) property of such approximations of singular
parabolic problems.






 As a simpler counterpart, it follows that for $c \le
c_{\rm H}$ the limit of $\{u_\e\}$ always exists,
 \beq
   \label{2m8}
    u_\e(x,t) \to \bar u(x,t) \quad \mbox{as} \quad \e \to 0,
    \eeq
for any data $u_0 \in L^2(\O)$, where $\bar u(\cdot,t) \in
L^2(\O)$ for all $t >0$. This means {\em existence} of a (unique)
solution.

 The nonexistence result (\ref{2m7})
actually means that, in the supercritical range $c>c_{\rm H}$,
 \beq
 \label{2m9}
 |\bar u(x,t)| =  \infty
 \eeq
 can be treated as a ``proper solution" of the original problem
 (\ref{2m1})--({\ref{2m3}) for any $L^2$ data satisfying
 some extra conditions, e.g.,
the regularity and positivity at the origin (or, more generally,
non-orthogonality to a positive lineal, see below). Then, as we
have mentioned at the beginning of Introduction, (\ref{2m9}) means
that the extended semigroup of such solutions obtained by regular
approximations is {\em discontinuous} at $t=0$ for sufficiently
arbitrary initial data $u_0$.
 Note that  the behaviour of the ``approximating" family $\{u_\e(x,t)\}$
as $\e \to 0$ to get (\ref{2m9}) can be extremely oscillatory. For
instance, for data $u_0(x)$
 that are
oscillatory at $x=0$, we present in Section \ref{Sect4} an example
 of radial solutions with {\em non-uniform oscillatory blow-up},
 where $u_\e(x,t)$ has both limits $\pm \infty$ along
some subsequences, and the same holds even at the origin $x=0$.

\ssk

 In the present paper, we study the
case where the singularity of the potential is concentrated at the
single internal point $x=0 \in \O$, and this indeed affects our
final nonexistence conclusions. We do not consider the case of
singularities on the boundary, which can lead to other
regularization spectral properties and hence different and more
difficult nonexistence criteria.
 In the elliptic case, there exist related extended Hardy's
inequalities for smooth domains $\Omega$ corresponding to
potentials with inverse square singularity on the boundary
$\partial \Omega$, which were introduced in \cite{BM}. See also
\cite{HHL}, where further extensions via perturbations of such
potentials were presented, and also \cite{Gh04, Gh06} for delicate
nonlinear counterparts.

\subsection{Layout of the paper}

Sections \ref{SectG}--\ref{SectExt} are occupied with various
aspects of nonexistence analysis of parabolic models such as
(\ref{2m1}) and their time-dependent and nonlinear extensions. To
show further application, in Section \ref{SectNN}, we discuss some
rather simple  corollaries of our analysis being applied to
Schr\"odinger and hyperbolic PDEs
 such as
 $$
  \mbox{$
 {\rm i} \, u_t=-(-\D)^m u + \frac c{|x|^{2m}}\, u \quad \mbox{and}
 \quad  u_{tt}=-(-\D)^m u + \frac c{|x|^{2m}}\, u,
 $}
 $$
 together with their semilinear counterparts.

\section{General nonexistence theorem}
\label{SectG}

\subsection{Divergence theorem for a general potential}

In this section, we treat the nonexistence in a more general
setting. We consider the Cauchy--Dirichlet  problem for the
poly-harmonic equation with a general regularized potential,
 \beq
 \label{2m1G}
  \mbox{$
u_\e: \quad u_t ={\bf B}_\varepsilon u \equiv -(-\D)^m u + \frac 1
 {\varepsilon^{2m}} \,q\big(\frac x\varepsilon\big)u \quad \mbox{in}
 \quad \O \times \re_+ \quad
  (m \ge 1).
  $}
 \eeq
 Here the potential
 $q_\e(x) \equiv \frac 1
 {\varepsilon^{2m}} \,q\big(\frac x\varepsilon\big)$
  depends on the parameter $\e>0$, for instance, in a manner
 similar to that in (\ref{2m5}),
where $q(y)$ is now an arbitrary smooth function  in
$\mathbb{R}^N$ decaying at infinity. Therefore, (\ref{2m1G}) plays
a role of a regular approximation of the parabolic equation with
the {\em singular} potential $V(x)$ ($V(0)=\infty$ and
$|V(x)|<\infty$ for $x \not = 0$) such that, uniformly on compact
subsets in $\ren \setminus\{0\}$ and sufficiently fast,
 \beq
  \label{ss11}
   \mbox{$
   q_\e(x)
   =
\frac 1
 {\varepsilon^{2m}} \,q\big(\frac x\varepsilon\big)
  \to V(x) \quad \mbox{as} \quad \e \to 0^+
  \quad \big(\mbox{e.g.,} \,\,\, V(x)= \frac c{|x|^{2m}} \,\,\, \mbox{as in (\ref{2m1})}\big).
 $}
  \eeq
 We assume for $u_\e(x,t)$  zero Dirichlet conditions on the smooth boundary,
 \beq
 \label{2m2G}
  \mbox{$
 u_\e= \frac{{\partial}u_\e}{\partial \nu}=...=\frac{{\partial}^{m-1}u_\e}{\partial
 \nu^{m-1}}=0 \quad \mbox{on} \quad \partial \O \times \re_+,
 $}
  \eeq
  where $\nu$ is the unit outward normal to $\partial \O$,
   and prescribe initial data
    \beq
    \label{2m3G}
    u_\e(x,0)=u_{0\e}(x) \in L^2(\O)
    \eeq
    that, in general,  also depend on $\e$.

    \ssk

In what follows, $C$ and $c_\e$ denote various positive constants
which exact values are of no importance. Here,  $C$ is independent
of $\e$, while constants $c_\e$ have at most rational dependence
on $\e$, i.e.,
 $$
\e^P<c_\e<\e^{-P}, $$
 for some positive $P$ which can be arbitrarily large.
 In this calculus, we easily write $c_\e c_\e \ge c_\e$.
Such constants are not of importance and are negligible while
dealing with exponential factors such as ${\mathrm e}^{ 1/ \e}$,
${\mathrm e}^{ 1/ {\e^2}}$, or ${\mathrm e}^{ 1/ {\e^{2m}}}$, to
be treated more carefully.

First of all, using the above calculus, since the equation is
linear, we may always assume
 \beq
 \label{ee12}
c_\e \le \|u_{0\e}\|_{L^2} \le c_\e \quad \mbox{for \, $\e>0$ \,
small},
  \eeq
i.e., $u_{0\e}$ does not get large as $\e \to 0$,
 which is a
natural and not restrictive assumption.

\ssk

    Our main hypothesis on the potential is as follows:

    \ssk

    \noi{\bf Hypothesis (P):} {\em The limiting operator}:
    \beq
     \label{oplim}
     {\bf  B}  \equiv
    -(-\D_y)^m  + q(y)I  \quad \mbox{\em in}\ \,\,  \mathbb{R}^N
    \eeq
   {\bf has $M \ge 1$ positive eigenvalues with exponentially decaying eigenfunctions}, {\em and,
    for some its  eigenfunction $U_n(y)$,}
    \beq
    \label{2m6G}
     \mbox{$
    |\langle u_{0\e}(x),U_n(\frac x \e)\rangle| \ge c_\e
    \quad \mbox{\em for all small} \quad \e>0.
     $}
     \eeq

\ssk

  {\em In the following  particular cases, $(\ref{2m6G})$ is replaced by}:

  \ssk

    (i) {\em if $u_{0\e}(\e x) \to \d_0>0$ as $\e \to 0$ $($e.g., $u_0$ is
     independent of $\e$, $u_0(x)$ is continuous at $x=0$ and
    $u_0(0)=\d_0>0)$,   there exists
    $U_n(y)$  that}
     {\bf has nonzero mean in $\ren$}; and

\ssk
     (ii) {\em if $u_{0\e}(x)$ is supported in some ball $B_{\bar
     c_\e}(0)$,
     with $\frac{\bar c_\e}\e \to 0$
     and $\int u_{0\e}(x) \ge c_\e$,
       there exists an eigenfunction  $U_n(y)$ of
    $(\ref{oplim})$
     that}
 {\bf does not vanish at the origin $y=0$.}

    \ssk


 Our main nonexistence result is a
corollary of the following estimate:

\begin{theorem}
 \label{Th.3}
  Let $m \ge 1$,  $(\ref{2m4})$,
     and Hypothesis {\bf (P)} hold.
   Then the sequence $\{u_\e(x,t)\}$ of classical solutions of the
  approximating problem $(\ref{2m1G})$--$(\ref{2m3G})$
 diverges in $L^2(\O)$ as $\e \to 0:$
   \beq
   \label{oc33}
    \|u_\e(x,t)\|_{L^2(\Omega)}\geq  {\mathrm e}^{\frac C {\e^{2m}}\,
    t} \to \iy
    \quad \mbox{for any fixed} \quad t>0.
    \eeq
\end{theorem}




\subsection{Proof of Theorem \ref{Th.3}}
 \label{Sect3}

Fix a sufficiently small $\e>0$. Consider the corresponding
eigenvalue problem
 \beq
 \label{3.1G}
 {\bf B}_\e \psi = \l\, \psi, \quad \psi \in H^{2m}(\O)\cap
 H^m_0(\O).
  \eeq
 By $\s({\bf B}_\e)= \{\l_j^\e, \, j \ge 0\}$ and $\{\psi_j^\e, \, j \ge 0\}$
we denote the corresponding spectrum and orthonormal, complete,
and closed in $L^2(\O)$ eigenfunctions subset, \cite{BS}. Clearly,
\beq
   \label{eq33}
    \mbox{$
    u_\e(x,t)= \sum_{(j \ge 0)} c_j^\e \, \psi_j^\e(x) \, {\mathrm e}^{\l_j^\e
    t}, \,\,\, \mbox{with} \,\,\, c_j^\e= \langle u_{0\e}, \psi_j^\e
    \rangle,
    $}
     \eeq
     \beq
     \label{eq33n}
 \mbox{$
   \mbox{so that} \quad \|u_\e(\cdot,t)\|_{L^2(\Omega)}^2=\sum_{(j \ge 0)} |c_j^\e|^2 \,
   {\mathrm e}^{2\l_j^\e
    t}.
  $}
    \eeq

We will need the following statement:

\begin{lemma}
\label{l33}
   Let $\Lambda_j,\ U_j$ for $j=1,2,...M$,  be positive eigenvalues
and  corresponding normalized eigenfunctions of
\eqref{oplim}. Then, for any small enough $\e>0$, intervals
 $$
  \mbox{$
 \big[ \frac {\Lambda_j}{\e^{2m}}-{\mathrm e}^
{-\frac C\e},\,  \frac {\Lambda_j}{\e^{2m}}+{\mathrm e}^ {-\frac
C\e}\big]
 $}
 $$
contain  eigenvalues of the problem \eqref{3.1G}. Moreover, the
function $U_j^\e=U_j \big(\frac x \e \big)$ is an ``approximate"
eigenfunction of $(\ref{3.1G})$ in the following sense:
  \be
   \label{GG1}
   \mbox{$
  \|U_j^\e-\sum_{(1)}\alpha_j\psi^\e_j\|_{H^m(\Omega)}\leq
  {\mathrm e}^
{-\frac C\e},
 $}
  \ee
where $\sum_{(1)}$ stands  for the summation with respect to $j$
such that $\lambda_j^\e$ belongs to the above interval, and
$\sum_{(1)}\alpha_j^2=c_\e$, $\alpha_j$ being constants
$($depending on $\e)$.
\end{lemma}

Let us return to the proof of the theorem. Due to the assumptions
of the theorem, we fix the necessary eigenfunction $U_n$  of the
operator \eqref{oplim} that provides us with the
  estimate
 \beq
  \label{GG2}
|\langle u_{0\e},U_n \rangle|\geq
c_\e.
 \eeq
  On
the other hand,  due to the above lemma, we see that there is at
least one eigenfunction $\psi^\e_l$ of \eqref{3.111} (see below)
with the corresponding eigenvalue $\lambda^\e_l$ in the interval
 $$
  \mbox{$
 \lambda^\e_l \in \big[\frac {\Lambda_j }{\e^{2m}}-C{\mathrm e}^
{-\frac C\e}, \frac {\Lambda_j }{\e^{2m}}+C{\mathrm e}^ {-\frac
C\e}\big].
 $}
 $$
In view of (\ref{GG1}) we have,
 \be
   \label{GG11}
    \mbox{$
   |\langle u_{0\e},U_n \rangle|-{\mathrm e}^ {-\frac {C}{\e}} \|u_{0\e} \|\leq
\big|\sum_{(1)}\alpha_l\langle u_{0\e},\psi^\e_l\rangle \big| \leq
\big(\sum_{(1)}\alpha_l^2\big)^{\frac 12}\big(\sum_{(1)}|\langle
u_{0\e},\psi^\e_l\rangle|^2\big)^{\frac 12},
 $}
  \ee
 and consequently, (\ref{ee12}) and \eqref{GG2} immediately implies that
  $$
  \mbox{$
\sum_{(1)}|c^\e_l|^2 \geq c_\e. $}
 $$
Then obviously we have
    $$
     \textstyle{
    \|u_\e(x,t)\|_{L_2(\Omega)}^2=
    \sum_{(j \ge 0)} |c_l^\e|^2 {\mathrm e}^{2\l_l^\e
    t}
    \geq\sum_{(1)} |c_l^\e|^2 {\mathrm e}^{2\l_l^\e t}
    \geq c_\e {\mathrm e}^{\frac{2\Lambda_j t}{\e^{2m}}}\sum_{(1)}
    |c_l^\e|^2\geq c_\e{\mathrm e}^{\frac{2C t}{\e^{2m}}}.
    }
    $$
  This
concludes the proof.
 $\qed$


\subsection{On generic asymptotic behaviour as $\e \to 0$}

 The following sharp
``pointwise" estimate of the solution sequence $\{u_\e\}$ holds:
if (\ref{2m6G}) is valid for $j=0$, then, for arbitrarily small
fixed $t>0$,
 \beq
 \label{as1n}
 \mbox{$
 {\mathrm e}^{-\l_0^\e \, t} u_\e(x,t) - c_0^\e \,
 \psi_0^\e(x) \to 0 \quad \mbox{as \,$\e \to 0$ \, in $L^2(\O)$}
 \quad \big( \,\psi_0^\e(x) \sim \e^{-\frac N2}\, U_0\big(\frac x \e \big)\, \big).
 $}
 \eeq




\subsection{A weaker blow-up hypothesis}
 The blow-up result
(\ref{oc33}) remains valid under the weaker condition on initial
data: (\ref{ee12}) holds and
 \beq
 \label{wc1}
  \mbox{$
 |\langle u_{0\e}(x), U_j \big( \frac x\e \big) \rangle | \ge  {\mathrm
 e}^{-\frac {c_*}{\e}}\ \ \text{for some positive constant}\ c_*<c_j,
  $}
  \eeq
where the constant   $c_j$ gives
 the rate of exponential decay of
the eigenfunction  $U_n$, i.e.,
 \be\label{GG3}
   \|U_j(y) {\mathrm e}^{c_j|y|}\|_{H^{m}(\mathbb{R}^N)} <C,
   \quad \mbox{with} \quad c_j>0.
    \ee
Assumption \eqref{wc1}  is obviously satisfied in the pioneering
original paper \cite{BG84} for $j=0$, see Proposition \ref{Pr.Sp}
(iii).

\subsection{Proof of Lemma $\ref{l33}$}

 Consider  eigenvalues and
eigenfunctions of operator \eqref{oplim} $\Lambda_j,\ U_j$,
$j=1,2,..., \, M$. These eigenfunctions decay exponentially at
infinity, see \eqref{GG3}.
     Then functions
$U_j(\frac x\e)$ satisfy the  equation
   \beq
 \label{3.111}
  \mbox{$
 {\bf B}_\e U_j\big(\frac x\e\big) = \frac 1{\e^{2m}}\, \Lambda_j U_j\big(\frac x\e\big)
  \quad \textrm{in} \,\,\,
  \Omega.
  $}
  \eeq
  These  do
not satisfy boundary conditions on $\partial\Omega$. To fix this,
consider the functions
 $$
  \mbox{$
 V_j(x)=U_j\big(\frac x\e \big)\chi(x),
  $}
  $$
   where $\chi$ is
smooth function which is equal to $1$ in a neighborhood of the
origin and is equal to zero in some neighborhood of the boundary
$\partial \O$. Then \beq \label{3.112}
 \mbox{$
 {\bf B}_\e V_j= \frac 1{\e^{2m}} \,\Lambda_j V_j+\tilde{V}_j \quad  \textrm{in} \
 \,\Omega, \quad
V_j \in H^{2m}(\O)\cap
 H^m_0(\O), \quad \mbox{and}
  $}
  \eeq
\be
 \label{GG5}
  \|\tilde V_j\|_{H^{m}(\O)}\leq {\mathrm e}^{-\frac C\e}.
   \ee
    Then
application of the  ``Lemma on approximate eigenfunction" (see
\cite{AV}) and observation that the difference between $U_j$ and
$V_j$ can be estimated via the right-hand side of \eqref{GG5}
deliver the result of the lemma. $\qed$

\ssk

 Lemma \ref{l33} can be strengthened (though, in this form, we
are not going to use it):

\begin{proposition}
 \label{Pr.Conv}
 Let $\Lambda_j$ and $U_j$ are eigenvalues
 and eigenfunctions of the problem \eqref{oplim}.
 Then eigenvalues and eigenfunctions of the eigenvalue problem
 $(\ref{3.1G})$satisfy:
  for any fixed $1\leq j \leq M$,
  as $\e \to
 0$,
   $$
  \mbox{$
 \big|\l_j^\e - \frac {\Lambda_j}{\e^{2m}}\, \big|\leq C{\mathrm e}^{-c|\Lambda_j|^{1/2m}\frac 1\e},
  \quad
  \|\psi^\e_j-\sum_{(\kappa_j)}\alpha_k^\e U_k(\frac x\e)\|_{L^2(\Omega)}
  \leq C {\mathrm e}^{-c|\Lambda_j|^{1/2m}\frac 1\e},
   $}
   $$
 where the summation in $\sum_{(\kappa_j)}$ takes place over $k$ such that $\Lambda_k=\Lambda_j$, and
$\sum_{(\kappa_j)}(\alpha_k^\e)^2= c_\e$.
 \end{proposition}




\section{Problem for (\ref{2m5}): scaling, limit problem, and
limit operator}
 \label{Sect01}

We now explain in greater  detail application of the above general
results to our main model equation (\ref{2m1}) with the
approximation (\ref{2m5}).

\subsection{Scaling}

Thus,  we perform in (\ref{2m5}) the natural scaling of the
independent variables
 \beq
 \label{2.1}
 x=\e \, y, \quad t= \e^{2m} \, \t.
  \eeq
  This reduces (\ref{2m5}) to the equation with the regular potential,
 \beq
 \label{2.2}
\mbox{$
 u_\t ={\bf B}_1 u \equiv -(-\D)^m u+ \frac c{1+|y|^{2m}}\, u \quad \mbox{in}
 \quad \O_\e \times \re_+,
  $}
 \eeq
 which is posed in the domain
  \beq
  \label{2.3}
   \mbox{$
  \O_\e= \big\{y \in \ren: \,\,\,   \e
    y \in \O\big\} \quad \big(\, u\big|_{\partial \O_\e}=0 \,\big).
 $}
  \eeq
The initial data are now
 \beq
 \label{2.4}
 u_{0\e}(x) \mapsto u_{0\e}(\e y) \quad \mbox{in}
 \quad \O_\e.
  \eeq

\subsection{Limit problem}

The above parabolic problems (\ref{2.2})--(\ref{2.4}) are posed in
a family of expanding domains $\{\O_\e\}$ that embrace the whole
$\ren$ as $\e \to 0$. Therefore, it is natural to consider the
{\em limit problem}, which is the Cauchy problem for
 \beq
 \label{2.5}
\mbox{$
 v_\t ={\bf B}_1 v \equiv -(-\D)^m v+ \frac c{1+|y|^{2m}}\, v \quad \mbox{in}
 \quad \ren \times \re_+.
  $}
 \eeq

\subsection{Spectral properties of limit operator}

The limit parabolic problem for (\ref{2.5}) introduces the {\em
limit operator}
 \beq
 \label{2.7}
  \mbox{$
 {\bf B}_1=-(-\D)^m + \frac c{1+|y|^{2m}}\, I \quad \mbox{in}
 \quad L^2(\ren), \quad \mbox{where} \quad c>c_{\rm H}(m).
  $}
  \eeq
 Obviously, (\ref{2.7}) is semibounded,
 \beq
 \label{2.8}
  {\bf B}_1 \le  c \, I  \quad \mbox{in}
  \quad C_0^\infty(\ren),
   \eeq
  so it admits  Friedrichs' self-adjoint extension (denoted again by ${\bf
  B}_1$)
   \cite{BS} with
  the domain $H^{2m}(\ren)$ in view of the embedding
 \cite[p.~54]{Maz}
   \beq
   \label{2.9}
   \mbox{$
   \int\limits_{\ren} \frac{w^2}{1+|y|^{2m}} \le C \int_{\ren}|D^m w|^2
   \quad \mbox{in} \quad H^m(\ren).
 $}
    \eeq
      For convenience, we  present a detailed description of
      the
  necessary  spectral properties:

\begin{proposition}
 \label{Pr.Sp}
{\rm (i)} The spectrum of the operator $(\ref{2.7})$ with domain
$H^{2m}(\ren)$ comprises the continuous and the discrete ones,
 \beq
 \label{2.10}
 \s({\bf B}_1) = \s_c ({\bf B}_1) \cup  \s_p({\bf B}_1)
 = \{ \l \le 0\} \cup \{\l=\Lambda_j>0, \, j=0,1,2,...\}, \quad
 \mbox{with}
 \eeq
  \beq
  \label{2.11}
  c>  \Lambda_0 \ge  \Lambda_1 \ge \Lambda_2 \ge ...  \ge  \Lambda_k \ge  \Lambda_{k+1} \ge ...
  \,0, \,
    \eeq
    where each eigenvalue repeated as many times as its finite
multiplicity.

{\rm (ii)} Each eigenfunction $U_j(y)$ is exponentially decaying
at infinity,
 \beq
 \label{222}
|U_j(y)| \le A_j {\mathrm e}^{-\a_m \Lambda_j^{1/2m} |y|} \quad
\mbox{as} \quad y \to \infty,
 \eeq
 where $A_j$ and  $\a_m$ are positive constants.

{\rm (iii)} The first eigenfunction $U_0(y)$ is radially
symmetric. For $m=1$, it is positive,
 \beq
 \label{2.22}
  U_0(y) > 0 \quad \mbox{in} \quad \ren \quad (m=1).
   \eeq
For any $m \ge 2$, $U_0(y)$ has infinitely many sign changes.

  \end{proposition}

\noi{\em Proof.}  {\bf Part (i).} \underline{\sc Continuous
spectrum.} Fix some $\Lambda = - \mu^{2m}<0$ and consider
 \beq
 \label{b1}
 {\bf B}_1 U = \Lambda U= - \mu^{2m} U \quad \mbox{in}
 \quad \ren.
  \eeq
  Using Fourier Transform, it is not difficult to see that
  (\ref{b1}) has a radial solution with the following asymptotic
  behaviour: as $|y| \to \infty$,
  \beq
  \mbox{$
  U_\Lambda(y) \sim |y|^{-\frac{N-2m}2} \cos(\mu |y|+a_0) \quad (a_0
  \in \re).
  $}\label{as1}
  \eeq
For $m=1$, $U_\Lambda(y)$ is given by Bessel's function $J_\nu(\mu
|y|)$, with $\nu= \frac{N-2}2$, \cite[\S~23]{Vlad72}. Here,
$U_\Lambda \not \in H^{2m}(\ren)$ (and $\not \in H^m(\ren)$, $\not
\in L^2(\ren)$). Hence, $({\bf B}_1-\Lambda I)^{-1}$ is not
bounded for any $\l <0$, from whence $0 \in \s_c({\bf B}_1)$ by
closure.


\smallskip

 \underline{\sc Discrete spectrum.}
 We first discuss a general representation of all the eigenfunctions,
 which we will use later on. To this end, we  introduce the
  polar coordinates $x=(r,\s)$ in $\ren$, so that
 \beq
 \label{L-B}
  \mbox{$
 \D= \D_r + \mbox{$\frac 1{r^2}$} \, \D_\s,
 \quad \mbox{where} \quad \D_r= D^2_r + \frac{N-1}r \, D_r
  $}
 \eeq
and $\D_\s$ is the Laplace--Beltrami operator on the unit sphere
$S^{N-1}$ in $\ren$, which is a regular operator with discrete
spectrum in $L^2(S^{N-1})$ (each eigenvalue repeated as many times
as its multiplicity)
 \beq
 \label{LBsp}
 \s(\D_\s) = \{- \mu_k \equiv -k(k+N-2) \le 0, \,\, k \ge 0\}.
 \eeq
 $\D_\s$ has  an orthonormal, complete, and closed subset $\{f_k(\s)\}$ of eigenfunctions being
$k$-th order homogeneous harmonic
 polynomials restricted to $S^{N-1}$.

Consider the eigenvalue problem for the operator (\ref{2.2}),
 \beq
 \label{l1}
   \begin{matrix}
 {\bf B}_1 U=(-1)^{m+1}\big(\D_r + \mbox{$\frac 1{r^2}$} \, \D_\s\big)^m  U +
 \frac c{1+|y|^{2m}}
 \,  U \smallskip\smallskip\\
 \equiv (-1)^{m+1} \sum_{l=0}^m C_m^l \D_r^l \big(\frac 1{r^2} \,
 \D_\s\big)^{m-l} U + \frac c{1+r^{2m}}
 \, U
  = \Lambda \, U,
  \end{matrix}
  \eeq
  where $C_m^l$ are binomial coefficients.
  Performing in (\ref{l1}) separation of variables
   \beq
   \label{l2}
   U(r,\s) = \phi(r) f_k(\s)
 \eeq
 yields the purely radial eigenvalue problem for $\phi(r)$,
  \beq
  \label{l3}
   \mbox{$
(-1)^{m+1} \sum_{l=0}^m C_m^l (-\mu_k)^{m-l} \D_r^l \big(\frac
1{r^{2(m-l)}} \,\phi
 \big) + \frac c{1+r^{2m}}
 \,  \phi
  = \Lambda \phi.
 $}
 \eeq
 Notice that the first operator in (\ref{l3})
 contains singular potentials up to the leading singularity in
 the last term
  $$
  \mbox{$
...  - \frac{c_{2m}}{r^{2m}} \, \phi, \quad \mbox{with a constant}
  \quad c_{2m}=\mu_k^m>0.
   $}
 $$
 This and other singular terms have the right sign to guarantee that this
 operator in negative and coercive, and, of course, this essentially affects the final
 result.

\ssk


 \underline{\sc
The positive lineal ${\mathcal L}_+= {\rm Span} \,\{U_0,U_1,...\}$
is infinite-dimensional.} The rest of the results are standard in
elliptic theory; see \cite{BS}. In particular, the existence of a
countable set of eigenvalues follows from the fact that the
subspace on which the quadratic form  for the operator
(\ref{2.7}),
 \beq
 \label{QQ1}
  {\bf Q}_1(v)= \langle v, {\bf B}_1 v \rangle,
   \eeq
is positive, is infinite-dimensional. Indeed, let $\chi(t)\in
C_0^\infty(\mathbb{R})$ be the cut-off function,
 $$
\chi(t)= \left\{
 \begin{matrix}
1, \,\,\, 0<t<1, \\
 0, \,\,\, t > 2, \,\,\,\,\,\,\,\,\,\,
  \end{matrix}
  \right.
  $$
and set
 \beq
 \label{CH1}
 \chi_{a,b}(t)=\chi(t-b)[1-\chi(t-a)], \,\,\, \mbox{$0<a<b$}
  \,\,\Longrightarrow \,\,
\text{supp} \ \chi_{a,b}\in [a,b+2].
 \eeq

 We next chose a test sequence for
Rayleigh quotion \cite{BS} in the following way:
  $$
u_{a,b}(|x|)=|x|^{-\frac{N-2m}{2}}\chi_{a,b}(|x|).
 $$
 Then we
have
 \beq
 \label{QQ4}
  \begin{matrix}
  {\bf Q}_1(u_{a,b})= \langle u_{a,b}, {\bf B}_1 u_{a,b} \rangle
  \ssk\ssk\ssk \\
   =\big\langle u_{a,b}, \big( \frac c{1+|x|^{2m}}- \frac{c_{\rm H}(m)}{|x|^{2m}}\big) u_{a,b}
  \big\rangle+\big\langle u_{a,b}, [-(-\D)^m,\chi_{a,b}] |x|^{-\frac{N-2m}{2}}
  \big\rangle,
   \end{matrix}
   \eeq
   where $[\cdot,\cdot]$ is the commutator.
   Here,  we have used the fact that ({\em q.v.} Section \ref{Sect4})
 $$
 \mbox{$
 -(-\D)^m |x|^{- \frac{N-2m}2}= \frac{c_{\rm H}}{|x|^{2m}}\, |x|^{-
 \frac{N-2m}2},
 $}
 $$
 so that the function $[-(-\D)^m,\chi_{a,b}]$ in (\ref{QQ4}) has
 support restricted to outer and inner layers according to (\ref{CH1}).
Obviously, the first term on the right-hand side in \eqref{QQ4}
tends to $+\infty$ as $b\rightarrow+\infty$ since $c>c_{\rm H}(m)$
and
 $$
 \mbox{$
 \int_{|x|>a+1}|x|^{-2m}|x|^{-2\frac{N-2m}{2}}\, {\mathrm d}x
 \sim \int^\iy r^{N-1} r^{-N} \, {\mathrm d}r =+\infty.
 $}
 $$
  In contrast, the  second term  on the right-hand side in \eqref{QQ4} remains
bounded, since
 $$
 \text{supp}\
[-(-\D)^m,\chi_{a,b}]\subset \{a\leq|x|\leq
a+1\}\cup\{b+1\leq|x|\leq b+2\}.
 $$

  As the result, ${\bf
Q}_1(u_{a,b})\rightarrow+\infty$ as $b\rightarrow+\infty$, and
there is a $b(a)$ such that ${\bf Q}_1(u_{a,b(a)})>0$ and
$\text{supp}\ u_{a,b}\in\{a\leq|x|\leq b+2\}$. Now choosing a
sequence $\{a_i\}$ such that $a_i>b(a_{i-1})+2$ for $i=1,2,...\,$,
we obtain infinite sequence $\{u_{a_i,b(a_1)}\}$
 such
that ${\bf Q}_1(u_{a_i,b(a_i)})>0$ and $\text{supp}\
u_{a_i,b(a_i)}\cap\text{supp}\ u_{a_j,b(a_j)}=\emptyset $ for
$i\neq j$. So there is an infinite dimensional subspace  on which
the quadratic form  for the operator (\ref{2.7}) is positive, and
consequently we have infinitely many positive eigenvalues of
$\BB_1$.

\ssk

 {\bf Part (ii).} Exponential decay follows from separation of variables and ODE
techniques. {\bf Part (iii).} The positivity (\ref{2.22}) for
$m=1$ is Jentzsch's classic theorem (1912).
$\qed$

\subsection{Nonexistence theorem}

Thus, we consider the Cauchy--Dirichlet problem for (\ref{2m5})
with general  data (\ref{2m3G}). Without loss of generality, we
state the following nonexistence (divergence) result, where, for
convenience, we slightly change the argument of the proof. This
underlines extra features of the above {\em positive lineal} of
the limit operator (\ref{2.7}),
 \beq
 \label{ll1}
 {\mathcal L}_+= {\rm Span}\{U_0, U_1,...\},
  \eeq
  as a set of all finite linear combinations of the given eigenfunctions.

\begin{theorem}
 \label{Th.NN}
 Let
  \beq
  \label{jj1}
   u_{0\e}(\e y) \to v_0(y) \quad \mbox{as \, $\e \to 0$\, in\,
   $L^2_{\rm loc}(\ren)$},
    \eeq
     and let, in the metric of $L^2(\ren)$,
     \beq
     \label{jj2}
      v_0 \,\,\, \mbox{is not orthogonal to} \,\,\,
      {\mathcal L}_+.
       \eeq
 Then $(\ref{oc33})$ holds.
 \end{theorem}

The proof repeats the arguments of the proof of Theorem
\ref{Th.3}. In particular, assuming the continuity and strict
positivity at the origin of the data that are independent of $\e$,
 \beq
 \label{rr1}
 u_{0\e}(x) = u_0(x) \quad \mbox{in} \quad \O,
 \quad \mbox{and} \quad u_0 \in C(B_1), \,\,\, u_0(0)=\d_0>0,
  \eeq
 we obtain from (\ref{jj1})
 \beq
 \label{2.6}
 v_0(y) \equiv \d_0=u_0(0)>0.
  \eeq
Then the condition (\ref{jj2}) (more precisely, see (i) below) is
obviously valid for $m=1$ in view of the positivity (\ref{2.22}).
For $m>1$, checking (\ref{jj2}) is not straightforward since all
the eigenfunctions are of changing sign. Anyway, we expect that
(\ref{jj2}) is always valid for constant data (\ref{2.6}). It is
curious that, for  $m \ge 2$, this is an open problem.

\subsection{On special solutions with fixed nodal sets}

Using the separation (\ref{l1}), we look for special solutions
 \beq
 \label{yy1}
 u(x,t)= v(r,t) f_k(\s), \quad r= |x|,
 \eeq
 of the original singular equation (\ref{2m1}) with $\O=B_1$. These solutions
 have fixed (time-independent) nodal sets that actually changes
the Hardy constant since now $0$ effectively belongs to an
artificial  boundary, on which $u=0$. Indeed, formally
substituting (\ref{yy1}) into (\ref{2m1}) yields the following
radial equation for $v(r,t)$:
 \beq
 \label{yy2}
  \mbox{$
 v_t=
(-1)^{m+1} \sum_{l=0}^m C_m^l (-\mu_k)^{m-l} \D_r^l \big(\frac
1{r^{2(m-l)}\, }\, v
 \big) + \frac c{r^{2m}}
 \,  v.
 $}
 \eeq
It follows that the last most singular at $y=0$ term is now
 \beq
 \label{yy3}
 \mbox{$
...+ \frac {c-(\mu_k)^m}{r^{2m}}
 \,  v,
  $}
 \eeq
 so that the problem becomes subcritical provided that
  \beq
  \label{yy4}
c-(\mu_k)^m \le c_{\rm H}, \quad \mbox{i.e., for} \,\,\, c \le
c_{\rm H} +(\mu_k)^m.
 \eeq
 Thus, for any $c> c_{\rm H}$,
there exist global solutions (\ref{yy1}) of (\ref{2m1}) with  $k
\gg 1$.

\ssk

Note that this by no means undermines the phenomenon of the strong
instability in the singular equation (\ref{2m1}) in the
supercritical range (\ref{2m4}). Indeed,  it is not clear which
solutions (\ref{yy1}), (\ref{yy4}) do withstand the
$\e$-regularization without blow-up as $\e \to 0$. We expect that,
without the assumption $\O=B_1$ or without special symmetries of
$\O$ supporting the fixed nodal set of special solutions under
consideration, no solutions actually do.
  However, in the
maximal generality, the proof is not a part of our business here.
In Section \ref{SectBl}, we present the proof of blow-up for some
stationary data that are supposed to be the most resistive  to
regular $\e$-approximations,
but actually are not.

\section{Constructing  data with doubly oscillatory $\pm \infty$ limit}
\label{Sect4}

In this section, for the reason of performing  some rather
involved calculus, we take the unit ball, $\O=B_1$, and restrict
to radial solutions only, where $r = |x| \in [0,1)$.

\subsection{Singular stationary solution (SSS)}

Consider the radial stationary equation
 \beq
 \label{r1}
  \mbox{$
 {\bf B}_0 U \equiv -(-\D_r)^m U + \frac c{r^{2m}}\, U=0 \quad
 \mbox{for} \quad r \in (0,1).
  $}
  \eeq
  As usual for Euler's-type ODEs,
looking for solutions of (\ref{r1})
 \beq
 \label{psi111}
 U(r) = r^\g,
 \eeq
 we obtain the following characteristic equation for $\g \in {\mathbb C}$:
 \begin{equation}
 \label{4.5l}
   \begin{matrix}
 G(\g) \equiv  G_*(\g) + c=0, \quad \mbox{where}
 \smallskip\smallskip\ssk\\
  G_*(\g)=(-1)^{m+1}
 \mbox{$\prod\limits_{k=1}^{m}$}
 [\g-2(k-1)](\g+N-2k).
 \end{matrix}
 \end{equation}
Then (\ref{Best1N}) means that (see e.g., \cite{GalH2})
 \begin{equation}
 \label{l2s}
 \mbox{$
 c=c_{\rm H} = -G_*(- \frac{N-2m}2) > 0,
 $}
 \end{equation}
 so the best Hardy constant $c_{\rm H}$ is such that
 the function (\ref{psi111}) with the exponent
 \beq
  \label{gam0}
  \mbox{$
\g= \g_m =  - \frac{N-2m}2
  $}
 \eeq
  is the exact weak solution of the
homogeneous equation $\BB_0 \psi = 0$ in $\ren$.
 Thus, for $c=c_{\rm H}$,
 the characteristic equation has the double root (\ref{gam0}).
It generates two $L^2$-solutions,
 \beq
 \label{2sing}
 \bar U_m(r) = r^{-\frac{N-2m}2} \ln r \quad \mbox{and} \quad
 \bar U_{m+1}(r) = r^{ -\frac{N-2m}2},
 \eeq
 which are ordered relative to the growth rate as $r \to 0$.
Other $2m-2$ characteristic roots of (\ref{4.5l}) are real or
complex.
Complex roots can occur for $m \ge 3$. One can see from the
structure of the characteristic polynomial (\ref{4.5l}) that, for
$N \gg 2m$ and $m$ even, there exist precisely two more real roots
$ \hat \g_{m-1}
> 2(m-1)$ and $\g_{m-1}< 2-N$. On the other hand, for $N \gg 2m$
and $m$ odd no more real roots exist; see \cite{GalH2}.

 It follows from (\ref{4.5l}) that in the supercritical range $c>
c_{\rm H}$, the characteristic polynomial admits two complex
roots,
 \beq
 \label{k1}
  \mbox{$
 \g_{m\pm}= -\frac{N-2m} 2 \pm {\mathrm i} \, d,
 \quad \mbox{where}
 \quad d =O\big (\sqrt{c-c_{\rm H}}\, \big)>0.
  $}
  \eeq
The corresponding solutions are oscillatory near the origin, e.g.,
 \beq
 \label{k2}
U_m(r) = r^{-\frac{N-2m}2}\cos(d \ln r + a_0) \quad (a_0 \in \re).
 \eeq

By a {\em weak SSS}  of (\ref{r1}) denoted by $U_*(r)$, we mean a
solution that exhibits the oscillatory behaviour such as
(\ref{k2}) as $r \to 0$. Note that
 \beq
 \label{k3}
  \mbox{$
 U_* \in L^p(\O) \quad \mbox{for any}
  \quad p < \frac {2N}{N-2m}.
  $}
  \eeq

\subsection
{$\e$-approximation diverges to $\pm \infty$}
 \label{S**}

In the radial parabolic approximating  problem for equation
(\ref{2m5}), we first take initial data $u_0(r)$ that are not
continuous at $r=0$ and exhibit an oscillatory behaviour. It is
convenient to assume the behaviour  as in (\ref{k2}), i.e.,
 \beq
 \label{u11}
 u_0(r) =  r^{-\frac{N-2m}2}\cos(d \ln r) \quad \mbox{for $r>0$
 small} \quad (\mbox{or, simply, $u_0(r)=U_*(r)$}),
  \eeq
though other types of oscillations will also do.



 Let us explain the origin of such huge oscillatory properties of
 the approximation $u_\e(x,t)$.
  Using the same idea, we apply
    the estimates from Section \ref{SectG} that are sufficient.
     The simpler
    case $m=1$, where $U_0(y)>0$ ({\em q.v.} (\ref{2.22})) was studied  in \cite{GBG1}, so here $m
    \ge 2$.

Thus,
   we begin with  formal estimates of
   the first Fourier coefficient for the  maximal
 positive  eigenvalue as in (\ref{as1n}) (here $r=|x|>0$),
    \beq
    \label{ss1}
     \mbox{$
    \l_0^\e
   \approx
   \frac {\Lambda_0}{\e^{2m}}
 \quad \mbox{and} \quad \psi_0^\e(x)
 \sim \e^{-\frac N2} U_0\big( \frac x \e \big)
 \sim \e^{-\frac N2}{\mathrm
 e}^{-a_m\Lambda_0^{1/2m}\, \frac r \e} \, \cos\big(b_m \, \frac r \e\big)
  \,\,\, \mbox{for} \,\,\, r \gg \e,
  $}
  \eeq
  where $\mu_m=\frac 1{\e}\,(-a_m + {\rm i} \, b_m) $ is the root of
   $(-1)^{m+1} \mu^{2m}= \frac {\Lambda_0}{\e^{2m}}$ with the
  maximal ${\rm Re} \, \mu_m=-a_m<0$.
This  yields, in the first rough, but sufficient, approximation,
the coefficient
    \beq
 \label{3.511}
  \mbox{$
c_0^\e= \langle u_0, \psi_0^\e \rangle \sim \e^{-\frac N2}
\int\limits_0^{\frac 1\e} r^{N-1} u_0(r) U_0\big(\frac r \e\big)\,
 {\mathrm d}r=
\e^{\frac N2} \int\limits_0^{\frac 1{\e^2}} z^{N-1} u_0(\e z) U_0(z)\,
 {\mathrm d}z.
  $}
   \eeq
Substituting next  data (\ref{u11}) yields the following estimate
of the first Fourier coefficient:
   \beq
 \label{3.512}
  \mbox{$
c_0^\e \sim \e^{m} \int\limits_0^{+\infty} z^{\frac{N+2m-2}2} U_0(z)\,
\cos(d \ln z - D) \,
 {\mathrm d}z \equiv \e^m I_d(D),
  $}
   \eeq
   where $ D=-d \ln \e \to +
 \infty$ as $\e \to 0$ is a shifting argument in $\cos(\cdot)$, which  changes the
 sign of the function in the integral.
 Recall that $U_0(z)$ is exponentially small by (\ref{ss1}) as $z
 \to \iy$, so the integral in (\ref{3.512}) fast converges.
 Here, we  assume that both integrals
  $$
  \textstyle{
 \int\limits_0^{+\infty} z^{\frac{N+2m-2}2} U_0(z)\,
\cos(d \ln z) \,
 {\mathrm d}z \quad \mbox{and}
 \quad  \int\limits_0^{+\infty} z^{\frac{N+2m-2}2} U_0(z)\,
\sin(d \ln z) \,
 {\mathrm d}z
 }
 $$
 do not vanish simultaneously (otherwise, we change $d$ and/or the
 exponent $-\frac{N-2m}2$ or any other parameters to get
 necessary non-zero values).
 Then,
in view of
 analytic dependence on parameters in (\ref{3.512}), for
any fixed $d>0$, there exists sequences $\{\e_k^{\pm}\} \to 0$
such that $c_0^{\e_k^+} \ge (\e_k^+)^m \d_0$ and $c_0^{\e_k^-} \le
-({\e_k^-})^m \d_0$ for any $k=0,1,2,...$ with some constant
$\d_0>0$. Therefore,
 by (\ref{ss1}), for $\e=\e_k^\pm \to 0$,
 \beq
 \label{91}
  \mbox{$
 u_\e(x,t_0) \sim c_0^\e\, \e^{m- \frac N2} I_d(d \ln \e)\,
{\mathrm
 e}^{-a_m\Lambda_0^{1/2m}\, \frac r \e} \, \cos\big(b_m \, \frac r \e\big)
 \,\, {\mathrm e}^{{
 \Lambda_0}\,\frac{t_0}{\e^{2m}}} \quad (r=|x|),
  $}
  \eeq
  where we may assume that $\cos(\cdot)$ does not violate the
  above sigh restrictions on $c_0^\e$.
 Similarly, in view of the oscillatory functions $\cos\big(b_m \, \frac r
 \e\big)$ in (\ref{91}), there exist sequences $\{\hat \e_k^{\pm}\} \to 0$
 such that $\psi_0^\e(x)$ stands positive and negative
 respectively.

Thus, we observe a {\em doubly oscillatory} singular limit to $\pm
\iy$ as $\e \to 0$ for the regularizing solution (\ref{91}).
 Imposing necessary elementary assumptions on these four sequences
  $\{\e_k^{\pm}\}$ and $\{\hat \e_k^{\pm}\}$,
 one always can choose some ``intermediate" sequences $\{\bar \e_k^{\pm}\} \to 0$
 along which the limits $\pm \iy$ are guaranteed.
     Thus, by (\ref{91}),
 we observe that, for any $r=|x|>0$,
  in the sense of the first Fourier coefficient,
 \beq
 \label{001}
 \mbox{$
 \limsup_{\e \to 0} \, u_\e(x,t) = + \infty
  \quad \mbox{and}
  \quad
 \liminf_{\e \to 0} \, u_\e(x,t) = - \infty
 $}
   \eeq
In particular,
 we also easily fix the following weaker
$\pm\iy$-divergence:
 \beq
 \label{001W}
 \mbox{$
 \limsup_{\e \to 0} \, \sup_x \, u_\e(x,t) = + \infty
  \quad \mbox{and}
  \quad
 \liminf_{\e \to 0} \, \inf_x\, u_\e(x,t) = - \infty
 $}
   \eeq
It is interesting
 that, by oscillatory behaviour of
the
 Fourier coefficient $c_0^\e$ in (\ref{3.512}), the non-uniform
divergence (\ref{001}) holds even at the origin $x=0$.

 On the
other hand, if for  given data $u_0$, the coefficient $c_0^\e$ is
not oscillatory at all, e.g.,
 $$
 c_0^\e  \ge c_\e>0 \quad \mbox{for all small} \quad \e>0,
  $$
the oscillatory divergence (\ref{001W}) or (\ref{001})  for any $x
\not = 0$ holds due to the sign-changing behaviour of $U_0(y)$ for
$y \gg 1$ given in (\ref{ss1}).







\section{An analytic stationary profile blows up as $\e \to 0$}
\label{SectBl}

Let us present a different simple example showing another
unstability feature of the $\e$-approximation of (\ref{2m1}),
which we now consider in $\ren \times (0,T)$, i.e., for
convenience, we pose the Cauchy problem.
 We  choose
$c>c_{\rm H}$ such that the characteristic equation (\ref{4.5l})
has the root $\g=2m$, so that (\ref{r1}) has the analytic
stationary solution
 \beq
 \label{ss1NN}
 U_0(x) = |x|^{2m} \quad \mbox{in} \quad \ren.
  \eeq

Consider the regularized equation (\ref{2m5}) with the stationary
data (\ref{ss1NN}). Choosing, for convenience, the new variable
 $$
 v=u_t
  $$
  and differentiating (\ref{2m5}) in $t$,
  we obtain for $v$ the same equation (\ref{2m5}) with  data
   \beq
   \label{ss2}
    \begin{matrix}
   v_{0\e}(x)= \BB_\e |x|^{2m} \equiv -(-\D)^m |x|^{2m}+ \frac
   c{\e^{2m}+|x|^{2m}}\, |x|^{2m}
   \ssk\\
   \equiv c\, \big( \frac 1{\e^{2m}+|x|^{2m}}- \frac
   1{|x|^{2m}}\big)=
    - \frac{c \,\e^{2m}}{\e^{2m}+|x|^{2m}}.
    \end{matrix}
 \eeq
 After scaling (\ref{2.1}), we arrive at the equation (\ref{2.5})
 with $\e$-independent initial data
  \beq
  \label{ss3}
   \mbox{$
  v_0(y)= -\frac c{1+|y|^2} < 0 \quad \mbox{in} \quad \ren.
   $}
   \eeq
Therefore, (\ref{jj2}) guarantees blow-up as $\e \to 0$. As usual,
 $m=1$  is easy by (\ref{2.22}).

Thus, in general, even classical stationary solutions of the
singular equation (\ref{2m1}) do not stand $\e$-regularization, to
say nothing about other more arbitrary (non-steady) data.

\section{Parabolic extensions: time-dependent and nonlinear potentials}
 \label{SectExt}

 \subsection{Time-dependent potentials: towards non self-adjoint
 spectral theory and a new Hardy constant}

 New spectral phenomena appear when the potential in (\ref{2m1G})
 depends on the time-variable in an $\e$-scaling manner, e.g.,
as in equation with the following singular {\em critical}
potential (cf. (\ref{2m1})):
 \beq
 \label{tt21}
  \mbox{$
 u_t= -(-\D)^m u + \frac {c}{t+|x|^{2m}}\, u.
  $}
  \eeq
  Then the natural regularized equation reads
 \beq
 \label{tt22}
  \mbox{$
u_\e: \quad  u_t= -(-\D)^m u + \frac {c}{\e^{2m}+t+|x|^{2m}}\, u.
  $}
  \eeq

As usual, we treat a general regularized equation with such
potentials,
  \beq
  \label{E.1}
   \mbox{$
u_\e: \quad  u_t= {\bf B}_\e u= -(-\D)^m u + \frac 1{\e^{2m}} \,
q\big(\frac x \e,\frac
  t{\e^{2m}}\big) \, u \quad \mbox{in}
  \quad \ren \times \re_+.
   $}
   \eeq
For convenience and as a different example, we  now ignore
boundary conditions and  consider the Cauchy problem with data
$u_{0\e} \in L^2(\ren)\cap L^\infty(\ren)$. Scaling (\ref{2.1})
yields for
 \beq
 \label{E.11}
  \mbox{$
 u_\e(x,t)= v_\e\big(\frac x \e, \frac t{\e^{2m}}\big)
  $}
 \eeq
 the following rescaled equation:
 \beq
 \label{E.2}
 v_\e: \quad
 v_\t= -(-\D)^m v + q(y,\t) v, \quad v(y,0)=v_{0\e}(y)=u_{0\e}(\e y).
  \eeq
Since by (\ref{E.11}) the behaviour of $\{u_\e\}$ as $\e \to 0$ is
equivalent to the behaviour of $\{v_\e(y,\t)\}$ as $\t \to +
\infty$, we perform  extra scaling and introduce the new standard
similarity variable from the heat kernel of the operator $D_t +
(-\D)^m$:
 \beq
 \label{E.3}
  \mbox{$
 z=\frac y{(1+\t)^{1/{2m}}}, \quad s= \ln (1+\t) \to +\infty.
  $}
  \eeq
Then
 $
 w_\e(z,s)= v_\e(y,\t)
 $
 solves the equation
  \beq
  \label{E.4}
  \mbox{$
  w_\e: \quad
  w_s = -(-\D)^m w + \frac 1{2m} \, z \cdot \n w +
   {\mathrm e}^s q \big(z {\mathrm e}^{\frac s{2m}}, {\mathrm e}^s-1 \big)\, w.
   $}
   \eeq

In order to approach a fixed regularized potential (as in
(\ref{2.5})), we assume that
 \beq
 \label{E.5}
  \mbox{$
  {\mathrm e}^s q \big(z {\mathrm e}^{\frac s{2m}}, {\mathrm e}^s-1 \big)
  \to \frac c{1+|z|^{2m}} \quad \mbox{as}
  \quad s \to + \infty \,\,\,\, \mbox{uniformly}
   $}
   \eeq
   and sufficiently 
   fast.
For instance, (\ref{E.5}) holds for the  potential in (\ref{tt22})
since
 \beq
 \label{E.6}
 \mbox{$
\frac 1{\e^{2m}} \, q\big(\frac x \e,\frac
  t{\e^{2m}}\big)= \frac c{\e^{2m}+t+|x|^{2m}}
  \quad \Longrightarrow \quad q(y,\t)= \frac c{1+\t+|y|^{2m}}.
   $}
   \eeq

Thus, passing to the limit in (\ref{E.4}) with the assumption
(\ref{E.5}), we obtain the equation
 \beq
 \label{E.7}
 \mbox{$
  w_s =
   -(-\D)^m w + \frac 1{2m} \, z \cdot \n w +
 \frac c{1+|z|^{2m}} \, w.
   $}
   \eeq
Here the limit operator
 \beq
 \label{E.8}
 \mbox{$
{\bf B}_1  \equiv -(-\D)^m  + \frac 1{2m} \, z \cdot \n  +
 \frac c{1+|z|^{2m}}\, I
 $}
 \eeq
admits  Friedrichs self-adjoint extension for the second-order
case $m=1$ only, where it possesses  a symmetric representation of
the form
 \beq
 \label{E.9}
\mbox{$ {\bf B}_1  \equiv \frac 1{\rho}\, \n \cdot (\rho \n) +
\frac c{1+|z|^{2m}}\, I
 \quad \mbox{in \,\, $L^2_\rho(\ren)$,\, with \,\, $\rho(z)={\mathrm
 e}^{\frac {|z|^2}4}$.}
 $}
 \eeq
For $m>1$, (\ref{E.8}) is not self-adjoint and is not symmetric in
 any weighted $L^2$-spaces.

Consider the linear operator with the same principal differential
part,
 \beq
 \label{bb1}
  \mbox{$
 {\bf B}=  -(-\D)^m  + \frac 1{2m} \, z \cdot \n  + \frac
 N{2m} \, I.
  $}
  \eeq
It is known \cite{Eg4} that, for any $m >1$, operator (\ref{bb1})
is naturally defined in the weighted space $L^2_\rho(\ren)$, with
 \beq
 \label{E.10}
  \mbox{$
\rho(y) ={\mathrm e}^{a|y|^\a}, \quad \mbox{with} \quad \a=
\frac{2m}{2m-1},
 $}
 \eeq
where $a>0$ is a sufficiently small constant. The domain of ${\bf
B}$ is the corresponding Sobolev space $H_\rho^{2m}(\ren)$, so
that $\BB: H_\rho^{2m}(\ren) \to L_\rho^{2m}(\ren)$ is bounded
with discrete spectrum
 \beq
 \label{bb2}
 \mbox{$
 \s({\bf B})= \big\{- \frac l{2m}, \,\,\, l=0,1,2,...\big\}.
  $}
   \eeq
   In particular,
 both ${\bf
B}$ and ${\bf B}^*$ then have discrete spectra, compact
resolvents, and complete sets of eigenfunctions. Some of these
results can be extended to the operator (\ref{E.8}), for which
 spectral theory deserves further
development.


\ssk

We now concentrate on particular spectral properties related to
the  nonexistence.

\noi\underline{\sc Consider first the self-adjoint case $m=1$}.
Then, by classic theory \cite{BS,Kato}, there exists a monotone
branch $\Lambda_0=\Lambda_0(c)$ of the first simple real
eigenvalue of $\BB_1$ (here $m=1$):
 \beq
 \label{mn1}
  \mbox{$
 \Lambda_0(0)= -\frac N{2m} \quad \mbox{and} \quad \Lambda_0(c) \to +
 \infty \,\,\, \mbox{as}
 \,\,\, c \to +\iy.
  $}
  \eeq
 We then  define the new ``Hardy
constant" for the $c$-family of operators (\ref{E.8}) as follows:
 \beq
 \label{bb3}
 c_{\rm H}: \quad \Lambda_0(c)=0 \in \s({\bf B}_1).
  \eeq
  Then, for any $c>c_{\rm H}$, there exists $\Lambda_0(c)>0$ that leads
to generic blow-up of approximations $\{u_\e\}$ provided that
$v_0$ is not orthogonal to $U_0$. Analogously, further increasing
$c>c_{\rm H}$ leads to appearance of more branches of positive
eigenvalues $\Lambda_1(c)$, ..., $\Lambda_K(c)$, so that
 we can have the positive lineal ${\mathcal L}_+=
{\rm Span}\,\{U_0, ...,U_K\}$ of arbitrarily large, but finite
dimension $K=K(c)$ for $c \gg 1$.

Thus, taking into account the first unstable mode only, we have
from scalings (\ref{E.11}) and (\ref{E.3}) that, under the
``non-orthogonality assumption" to $U_0 \in {\mathcal L}_+$, the
following divergence rate is achieved: for an arbitrarily small
fixed $t>0$ (here again $m=1$),
 \beq
 \label{ss1nn}
  \mbox{$
   u_\e(x,t) \sim \big( \frac t{\e^{2m}}\big)^{\Lambda_0(c)}
   U_0\big(\frac x{t^{1/2m}}\big) \,\,\, \mbox{as} \,\,\, \e \to
   0.
   $}
   \eeq
Similar to Theorem \ref{Th.3}, the behaviour (\ref{ss1nn}) is easy
to characterize in terms  of $L^p(\ren)$-divergence for any $p \ge
1$, as well as in other metrics.

 On the contrary,
 for $c \le c_{\rm H}$, the
family $\{u_\e=v_\e=w_\e\}$ is uniformly bounded meaning existence
of a bounded solution (of (\ref{tt21}), say).

\ssk

\noi\underline{\sc Consider  $m \ge 2$}, where we can guarantee
for sure less. Since in (\ref{E.8}) the last term $\frac
c{1+|z|^{2m}}\, I
$
 serves as a compact perturbation of $\BB$ in the equivalent
 integral representation, we still have a discrete spectrum and a countable
 family of continuous eigenvalue branches for all $c \ge 0$.
  For $c=0$, there holds:
 $$
  \mbox{$
 {\bf B}_1={\bf B}- \frac N{2m} \, I,
 $}
 $$
so by  perturbation methods \cite{Kato}, we find the local branch
as in (\ref{mn1}) of the first real {\em simple} eigenvalue
$\Lambda_0(c)$ of ${\bf B}_1$ for all sufficiently small $c>0$.
 Further extension of the branch in this non self-adjoint case is
 not guaranteed to be strictly monotone, and the eigenvalue does
 not necessarily remain
 always real. Then we are assumed to deal with the first
 eigenvalue having the maximal real part ({\em q.v.} (\ref{bb3})),
  \be
  \label{nat1}
  {\rm Re} \,
 \Lambda_0(c)>0 \,\,\, \mbox{for} \,\,\, c>c_{\rm H}.
  \ee
This definition is based on a natural property that ${\rm Re} \,
\Lambda_0(c)$ growths unboundently as $c \to + \iy$, though
passage to the limit $c \to +\iy$ in the family
 (\ref{E.8}) is quite tricky and not studied here.
 Then under the assumption (\ref{nat1}) (or (\ref{bb3}) in the real case), we still
obtain the divergence result of the type (\ref{ss1nn}) for  any
non-orthogonal initial data now in the sense that $\langle v_0,
U_0^* \rangle \not = 0$, where $U_0^*$ is the corresponding
 adjoint eigenfunction of $\BB^*_1$.


Thus, for all $m \ge 1$, according to our previous results, the
divergence of the family $\{u_\e\}$ as $\e \to 0$ will be
guaranteed under the hypothesis of the type (\ref{jj2}), where
${\mathcal L}_+$ is the positive lineal of the limit operator
(\ref{E.8}). Indeed, deeper
 spectral analysis of the non self-adjoint operators (\ref{E.8})
 are quite in demand here.

\subsection{On nonlinear parabolic equations}

As a next natural extension, we consider the following semilinear
parabolic equation with the same singular potential as in
(\ref{2m1}):
 \beq
 \label{eee1}
  \mbox{$
 u_t= -(-\D)^m u + \frac 1{|x|^{2m}} \, |u|^{p}, \quad
 \mbox{where}
 \quad p>1.
 $}
  \eeq
 Here, the nonlinearity is taken in the form $|u|^p$ for convenience of
further simple calculus in the case $m \ge 2$, when the Maximum
Principle fails.

The condition $p>1$ involves the extra high instability of the
evolution driving by the reaction-diffusion PDE (\ref{eee1}),
since now it described blow-up of solutions in finite time. For
$m=1$ and $u_0 \ge 0$, the nonexistence results for (\ref{eee1})
 were obtained in
Brezis--Cabr\'e \cite{BrCabre}; on full history and extensions
 see  \cite{GGK05, Kom06},
\cite[p.~333]{GalGeom}, and \cite[p.~108,\,267]{MitPoh}. For
solutions of changing sign, the nonexistence results for
(\ref{eee1}) were absent even for $m=1$.

 Let us
see how the nonlinearity $|u|^p$  affects the regularized
solutions satisfying, as usual,
 \beq
 \label{eee2}
  \mbox{$
   u_\e: \quad
 u_t= -(-\D)^m u + \frac 1{\e^{2m}+|x|^{2m}} \, |u|^{p}.
 $}
  \eeq
Then scaling (\ref{2.1}) yields
 \beq
 \label{eee3}
  \mbox{$
   v_\e: \quad
 v_t= -(-\D)^m v + \frac 1{1+|y|^{2m}} \, |v|^{p}
 \quad \mbox{in}
 \quad \ren \times \re_+,
 $}
  \eeq
and for simplicity we assume that data are constant as in
(\ref{2.6}).

\ssk

\noi\underline{\sc Case ${m=1}.$} Then the proof of blow-up of
(\ref{eee3}) is straightforward, since by (\ref{2.6}),
 \beq
 \label{vv4}
 v_t > 0  \quad \mbox{in}
 \quad \ren \times \re_+
 \eeq
 via the Maximum Principle. Therefore, the solution $v(y,\t)$ is
 strictly monotone increasing in $t$ and is not bounded above
 (since otherwise, as for gradient systems, $v(\t)$ would stabilize to
 a bounded steady solution,
 which is nonexistent); see details in \cite[p.~319]{GalGeom}.
  Denoting by $T_*$ the blow-up time for
 (\ref{eee3}), we take into account the fact that
due to (\ref{vv4}) blow-up is {\em complete}, i.e., in the natural
sense of the proper minimal extension,
 $$
 v(y,\t) \equiv +\infty \quad \mbox{for
\, $\t>T_*$}.
 $$
  We  then conclude that, for all small $\e>0$,
 \beq
 \label{eee4}
 u_\e(x,t) \,\,\, \mbox{blows up completely in time $T_\e <
 2\e^{2m}T_*$}.
  \eeq
This ends the proof of divergence as $\e \to 0$ of $\{u_\e\}$
under the assumption (\ref{rr1}) for $m=1$.

\ssk

\noi\underline{\sc Case ${m \ge 2}.$} A similar scheme applies
with some changes since (\ref{vv4}) is no longer true. Finite time
blow-up in (\ref{eee3}), (\ref{2.6}) is then proved by the
eigenfunction method that shows blow-up of the first Fourier
coefficient
 $$
  \mbox{$
 E(\t)= \big\langle v(y,\t),\psi_1\big(\frac y R\big)\big\rangle \quad (R \gg 1),
  $}
  $$
  where $\psi_1(z)$ is a specially constructed nonnegative cut-off function
  satisfying the following elliptic inequality: there exists a
  constant $\l_1>0$ such that
  \begin{equation}
\label{Ap2}
 |\Delta^m\psi_1| \le \lambda_1 \psi_1 \quad {\rm in}
\,\,\, \ren.
\end{equation}
See a convenient adaptation of Kaplan's eigenfunction method
(1963) to higher-order parabolic operators in \cite[\S~5]{GP1}.
 It is crucial that the
resulting ordinary differential inequality for $E(\t)$ implies not
only the fact of finite-time blow-up at some $T_*$ but also that
it becomes complete (possibly, at some moment $T_{\rm c} > T_*$).
 An alternative proof of blow-up in (\ref{eee3}) by the nonlinear
capacity method can be performed as in \cite[\S~29]{MitPoh}.
 This leads to the nonexistence conclusion (\ref{eee4}). Here,
 similar to (\ref{as1n}), we control {\em complete} blow-up
 of the ``first Fourier coefficient" relative to a solution
 $\psi_1$ of the elliptic inequality (\ref{Ap2}).


 The above scheme on nonexistence  applies
to the truly quasilinear parabolic model
 \beq
 \label{eee2n}
  \mbox{$
   u_\e: \quad
 u_t= -(-\D)^m (|u|^n u) + \frac 1{\e^{2m}+|x|^{2m}} \, |u|^{n+1}, \quad \mbox{where} \quad n>0.
 $}
  \eeq

\section{On nonexistence for  other linear and semilinear PDEs}
 \label{SectNN}

 \subsection{Schr\"odinger-type equations}

A typical linear example is
 \beq
 \label{2m5S}
  \mbox{$
  u_\e: \quad
 {\rm i} \, u_t ={\bf B}_\e u = -(-\D)^m u + \frac c{\e^{2m}+|x|^{2m}}\, u.
  $}
 \eeq
 The principle difference is that all $u_\e$ satisfy the
 conservation law
 \beq
 \label{cl1}
 \|u_\e(t)\|_{L^2}=\|u_{0\e}\|_{L^2} \quad \mbox{for all} \quad t >0.
  \eeq
Nevertheless, the nonexistence conclusion for $c>c_{\rm H}$ is
derived in a similar manner, though, in view of (\ref{cl1}), the
key feature of divergence is that the family $\{u_\e(x,t)\}$ gets
extremely oscillatory as $\e \to 0$ in the variable $t>0$. Indeed,
 now, loosely speaking, for positive eigenvalues $\Lambda_j$ of
${\bf B}_1$,
 \beq
 \label{tt1}
  \mbox{$
 \l_j^\e \sim -{\rm i}\, \frac{\Lambda_j}{\e^{2m}} \quad \mbox{for small} \quad \e>0.
  $}
  \eeq
Finally, instead of (\ref{as1n}), in view of the oscillatory
behaviour with purely imaginary eigenvalues (\ref{tt1}), assuming
 the strict inequality $\Lambda_0> \Lambda_1$ (cf. (\ref{2.11})),
we suggest  the following characterization of the actual
divergence of solution (of course this is just a formal
illustration that is not that informative, which any isolated mode
admits):
 \beq
 \label{as1nn}
 \mbox{$
 {\mathrm e}^{ {\rm i} \, \l_0^\e \, t} u_\e(x,t)
 -
  c_0^\e \,
 \psi_0^\e\big(\frac x \e\big) \rightharpoonup 0 \quad \mbox{as \,$\e \to 0$ \, weakly
  in $L^2_{\rm loc}(\O \times (0,\d))$.}
 $}
 \eeq
 Thus, in terms of $u_\e$, we get rather weak blow-up. However,
 the time-derivative $(u_\e)_t$ is then truly blows up as $\e \to 0$  according
 to
 $$
 \mbox{$
 (u_\e)_t(x,t) \sim \frac 1{\e^{2m}} \, {\mathrm e}^{ -{\rm i} \, \l_0^\e \, t}
 c_0^\e \,
 \psi_0^\e\big(\frac x \e\big).
  $}
  $$
The real parts of the solutions,
 $
 u_\e=U_\e+{\rm i}  V_\e,
 $
 satisfy a $4m$th-order hyperbolic equation,
  $$
  U_\e: \quad U_{tt}=- {\bf B}_\e^2 \, U,
  $$
  to which a general scheme  of the analysis, with some changes,
   can be also applied.




 As usual, a similar ``nonexistence" analysis of
  the corresponding {\em nonlinear Schr\"odinger equation} (NLS)
 \beq
 \label{2m5SN}
  \mbox{$
  u_\e: \quad
 {\rm i} \, u_t ={\bf B}_\e u \equiv -(-\D)^m u + \frac 1{\e^{2m}+|x|^{2m}}\, |u|^{p-1}u
 \quad (p>1)
  $}
 \eeq
 demands extra results on blow-up in the corresponding rescaled
 equation
\beq
 \label{2m5SN1}
  \mbox{$
  v_\e: \quad
 {\rm i} \, v_\t ={\bf B}_1 v \equiv -(-\D)^m v + \frac 1{1+|y|^{2m}}\,
 |v|^{p-1}v.
  $}
 \eeq
For $m=1$, blow-up for such  NLSEs in some parameter ranges of $p$
and $N$ is a classic matter of the theory and application to
nonlinear optics, plasma
 physics, and others; see  references to key papers and monographs on this subject
  in  Merle--Raphael \cite{Mer05} and in
  survey
\cite{GVSur02}. The higher-order case $m  \ge 2$
 is  less
developed though some techniques can be applied.
 Notice that we are interested in quite
 particular setting with constant data as in (\ref{2.6}),
 for which the proof of blow-up seems not  very
 difficult.

Thus, we will get typical ``divergence" of the family $\{u_\e\}$
as $\e \to 0$, provided that the suitable rescaled solution of
(\ref{2m5SN1}) blows up in finite time. Note that finite
$L^2$-energy blow-up solutions of NLSEs
 admit a natural
extension beyond blow-up time, i.e., for $t>T$ in view of the
conservation (\ref{cl1}) (possibly, this is not the case for data
(\ref{2.6}) of infinite energy). Then the divergence of $\{u_\e\}$
as $\e \to 0$ might be not that lethal as in other examples,
 so again a further study of such singularity phenomena is
necessary.



  \subsection{Hyperbolic equations}

  A key  regularized linear model is constructed analogously,
  \beq
  \label{hp1}
   \mbox{$
   u_\e: \quad
   u_{tt} ={\bf B}_\e u \equiv -(-\D)^m u + \frac c{\e^{2m}+|x|^{2m}}\,
   u,
 $}
  \eeq
with initial data $u_0$, $u_1$. Taking $t=\e^m \t$ in (\ref{2.1})
yields the rescaled equation
 \beq
  \label{hp2}
   \mbox{$
   v_\e: \quad
   v_{\t\t} ={\bf B}_1 v \equiv -(-\D)^m v + \frac c{1+|y|^{2m}}\,
   v.
 $}
  \eeq
The separation of variables
generates the eigenvalue problem (\ref{3.1G}), where $\l \mapsto
\l^2$.
 Hence,   the positive lineal of ${\bf
B}_1$ continues to play the key role for the nonexistence as $\e
\to 0$, so under similar non-orthogonality assumptions, we arrive
at the divergence result such as (\ref{as1n}), with
$\sqrt{\l_0^\e}$ in the first multiplier.

The corresponding nonlinear hyperbolic equation,
 \beq
  \label{hp1N}
   \mbox{$
   u_\e: \quad
   u_{tt} ={\bf B}_\e u \equiv -(-\D)^m u + \frac 1{\e^{2m}+|x|^{2m}}\,
   |u|^p,
 $}
  \eeq
leads to the rescaled PDE
 \beq
  \label{hp2N}
   \mbox{$
   v_\e: \quad
   v_{\t\t} ={\bf B}_1 v \equiv -(-\D)^m v + \frac 1{1+|y|^{2m}}\,
   |v|^p.
 $}
 \eeq
 For data (\ref{2.6}),
the necessary
  blow-up results for (\ref{hp2N}) can be
obtained as in  \cite[Ch.~5-7]{MitPoh}, so that the
nonexistence-divergence conclusions as $\e \to 0$ persist in a
similar way.





\end{document}